\documentclass[a4paper,12pt]{article}
\usepackage[utf8]{inputenc} 
\usepackage[english]{babel}
\usepackage{amsfonts,amssymb}
\title{Lyapunov-Krasovsky Functionals and the Existence of Unbounded Solutions of Functional-Differential Equations of Retarded Type \footnote{This research was supported by the Russian Science Foundation (project No. 25-21-00830)}}
\author{Lapin K.S., Lapin S.V.}
\date{}
\begin{document}
\maketitle

\hspace{11.7cm}UDC 517.925.54
\vspace{0.5cm}

\begin{abstract}
In this paper, the concepts of unboundedness and partial unboundedness of solutions of functional-differential equations of retarded type are introduced. On the basis of the method of Lyapunov-Krasovsky functionals, sufficient conditions for the existence of partially unbounded solutions are obtained. As corollaries, sufficient conditions for the existence of unbounded solutions are obtained.
\vspace{1mm}

Bibliography: 11 titles.
\vspace{1mm}


\end{abstract}
\vspace{0.5cm}

In the work of T. Yoshizawa \cite{Io1}, on the basis of the method of Lyapunov-Krasovsky functionals \cite{Kra}, a theory of boundedness of solutions of functional-differential equations of retarded type was created. This theory is an analogue of the classical theory of T. Yoshizawa of boundedness of solutions of systems of ordinary differential equations \cite{Io2} in the sense that in this theory, instead of a finite-dimensional space of initial positions of solutions, an infinite-dimensional functional space of prehistories of solutions is used. The main difficulties in the development of this theory of boundedness of solutions, as was shown in \cite{Io1} (see also \cite{H}), arise precisely in connection with the infinite-dimensionality of the functional space of prehistories of solutions. On the other hand, in the works of K.S. Lapin \cite{Lap1}-\cite{Lap5}, on the basis of the method of Lyapunov functions, sufficient conditions for the existence of various types of unbounded solutions of systems of ordinary differential equations were obtained. In this connection, an interesting and important problem for theory and applications arises: to obtain, on the basis of the method of Lyapunov-Krasovsky functionals, sufficient conditions for the existence of unbounded solutions of functional-differential equations of retarded type, which in the above-mentioned sense would be analogues of the sufficient conditions for the existence of unbounded solutions of systems of ordinary differential equations obtained in \cite{Lap1}-\cite{Lap5}. It is worth noting that at present there are very many works devoted to the study, using Lyapunov-Razumikhin functions, of various properties (stability, boundedness, etc.) of solutions of functional-differential equations of retarded type. However, starting from the work of T. Yoshizawa \cite{Io1} and up to the present time, no theoretical investigations using Lyapunov-Krasovsky functionals devoted to the unboundedness of solutions of functional-differential equations of retarded type have been carried out. Therefore, the relevance of solving the above-mentioned problem is also determined by the need to eliminate the existing gap connected with the use of Lyapunov-Krasovsky functionals for the study of unboundedness of solutions in the theory of functional-differential equations of retarded type.

In the present paper, for functional-differential equations of retarded type, the concepts of unbounded solutions and partially unbounded solutions are introduced. In terms of Lyapunov-Krasovsky functionals, sufficient conditions for the existence of partially unbounded solutions of functional-differential equations of retarded type are obtained (Theorems 1 and 2). As corollaries, sufficient conditions for the existence of unbounded solutions of functional-differential equations of retarded type are obtained (Corollaries 1 and 2). Then, in terms of two Lyapunov-Krasovsky functionals, sufficient conditions for the existence of partially unbounded solutions of functional-differential equations of retarded type are obtained (Theorems 3 and 4). As corollaries, sufficient conditions for the existence of unbounded solutions of functional-differential equations of retarded type are obtained (Corollaries 3 and 4). Further, sufficient conditions for the existence of partially finally-infinite solutions of functional-differential equations of retarded type are obtained (Propositions 1 and 2). As corollaries, sufficient conditions for the existence of finally-infinite solutions of functional-differential equations of retarded type are obtained (Corollaries 5 and 6). In conclusion, an example of the application of Theorem 1 to a concrete integro-differential system with positive constant time delay is given. Let us now turn to precise definitions and statements.

First, we recall the necessary concepts and constructions (see, for example, \cite{Io1}, \cite{H}) that are required for the consideration of functional-differential equations of retarded type.

Let $h>0$ be a fixed constant, and let $C([-h;0],\mathbb{R}^n)$ be the linear space of all continuous mappings of the form $\varphi:[-h;0]\to\mathbb{R}^n$, $n\geqslant 1$, considered together with the norm
$$\|\varphi\|=\max_{-h\leqslant\vartheta\leqslant 0}\|\varphi(\vartheta)\|,$$
where $\|\varphi(\vartheta)\|$ is the Euclidean norm in $\mathbb{R}^n$ of the vector $\varphi(\vartheta)\in\mathbb{R}^n$.
It is easy to see that any continuous mapping $x:[-h;A)\to\mathbb{R}^n$, $x(t)\in\mathbb{R}^n$, $-h\leqslant t< A$, where $A>0$ or $A=+\infty$, determines a continuous mapping $$\overline{x}:[0;A)\to C([-h;0],\mathbb{R}^n),\quad \overline{x}(t)=x_t\in C([-h;0],\mathbb{R}^n),\quad 0\leqslant t<A,$$ where the mapping $x_t:[-h;0]\to\mathbb{R}^n$ is defined by the formula
$x_t(\vartheta)=x(t+\vartheta)$, $-h\leqslant \vartheta\leqslant 0$.
It is clear that $x_t(0)=x(t)$ and $x_t(-h)=x(t-h)$.

Consider the functional-differential equation of retarded type
$$\dot{x}(t)=F(t,x_t),\quad t\in\mathbb{R}^+\equiv [0;+\infty),\eqno(1)$$
where $x(t)\in\mathbb{R}^n$, $-h\leqslant t<+\infty$, $F:\mathbb{R}^+\times C([-h;0],\mathbb{R}^n)\to \mathbb{R}^n$ is any continuous mapping, $\dot{x}(t)$ is the right-hand derivative of the mapping $x(t)$ and the family of mappings $x_t:[-h;0]\to \mathbb{R}^n$, $t\in\mathbb{R}^+$, is determined in the above-mentioned way from the mapping $x(t)$, $-h\leqslant t<+\infty$. It is worth noting that the mapping $x(t)$, $-h\leqslant t<0$, whose derivative is not present in the left-hand side of equation $(1)$, is required for the determination of the family of mappings $x_t:[-h;0]\to \mathbb{R}^n$, $0\leqslant t<h$, present in the right-hand side of equation $(1)$. For $n\geqslant 2$, equation $(1)$ is also called a system of functional-differential equations of retarded type.

Recall that a solution $x(t,t_0,\varphi_0)$ of equation $(1)$ with initial conditions $(t_0,\varphi_0)\in [0,A)\times
C([-h;0],\mathbb{R}^n)$, where $A>0$ or $A=+\infty$, is a mapping
$$x(-,t_0,\varphi_0):[t_0-h;A)\to\mathbb{R}^n,\quad x(t,t_0,\varphi_0)\in\mathbb{R}^n,\quad t_0-h\leqslant t<A,$$
for which the following conditions are satisfied:

1) $x(t,t_0,\varphi_0)$, $t_0\leqslant t<A$ satisfies equation $(1)$, i.e., for $t_0\leqslant t<A$ the equality $$\dot{x}(t,t_0,\varphi_0)=F(t,x_t(t_0,\varphi_0)),\quad\mbox{where}~~x_t(t_0,\varphi_0)=\overline{x(-,t_0,\varphi_0)}(t)=x(-,t_0,\varphi_0)_t,$$
holds;

2) $x_{t_0}(t_0,\varphi_0)=\varphi_0$, i.e., the equalities
$$x_{t_0}(t_0,\varphi_0)(\vartheta)=x(t_0+\vartheta,t_0,\varphi_0)=\varphi_0(\vartheta),\quad -h\leqslant \vartheta\leqslant 0,$$
hold.
The mapping $\varphi_0:[-h;0]\to\mathbb{R}^n$ is called the prehistory of the solution $x(t,t_0,\varphi_0)$ from the time moment $t_0-h$ to the time moment $t_0$ or, more briefly, the prehistory of this solution.

In what follows, we shall always assume that the continuous mapping $F(t,\varphi)$, which determines the right-hand side of equation $(1)$, satisfies a local Lipschitz condition with respect to the variable $\varphi$. The continuity of the mapping $F(t,\varphi)$ and the fulfillment of the indicated Lipschitz condition for this mapping ensure \cite{H}, respectively, the existence and uniqueness of the solution $x(t,t_0,\varphi_0)$ of equation $(1)$ with any given initial conditions $(t_0,\varphi_0)\in \mathbb{R}^+\times C([-h;0],\mathbb{R}^n)$.

Numerous examples of functional-differential equations of the form $(1)$ are given by integro-differential equations with positive constant time delay. Indeed, consider the integro-differential equation
$$\dot{x}(t)=f\left(t,x(t),x(t-h),\int_{t-h}^tg(s-t,x(s))ds\right),\quad t\in\mathbb{R}^+,\eqno(2)$$
with constant delay $h>0$, where $x(t)\in\mathbb{R}^n$, $-h\leqslant t<+\infty$, and the mappings $f:\mathbb{R}^+\times\mathbb{R}^n\times\mathbb{R}^n\to \mathbb{R}^n$ and $g:\mathbb{R}^+\times\mathbb{R}^n\to \mathbb{R}^n$ are any continuous mappings. For $n\geqslant 2$, equation $(2)$ is also called a system of integro-differential equations with constant delay. It is easy to see that every integro-differential equation $(2)$ with constant delay $h>0$ determines a functional-differential equation of retarded type $\dot{x}(t)=F(t,x_t)$,
$t\in\mathbb{R}^+$, where the continuous mapping $F:\mathbb{R}^+\times C([-h;0],\mathbb{R}^n)\to \mathbb{R}^n$ is defined by the following formula:
$$F(t,\varphi)=f\left(t,\varphi(0),\varphi(-h),\int_{-h}^0g(\vartheta,\varphi(\vartheta))d\vartheta\right).$$

Recall (see, for example, \cite{R}) that the upper right derivative of a continuous function $f:[a;+\infty)\to\mathbb{R}$ is the mapping $f'^+:[a;+\infty)\to\mathbb{R}\cup\{\pm\infty\}$ defined by the following formula:
$$f'^+(t)=\overline{\lim_{\delta\to\, 0^+}}\,\,\frac{1}{\delta}\,(f(t+\delta)-f(t)).$$
It is well known \cite{R} that a continuous function $f:[a;+\infty)\to\mathbb{R}$ is non-increasing if and only if $f'^+(t)\leqslant 0$ for all $t\geqslant a$. It is also well known \cite{Mc} that a continuous function $f:[a;+\infty)\to\mathbb{R}$ is non-decreasing if and only if $f'^+(t)\geqslant 0$ for all $t\geqslant a$. It is clear that if a continuous function $f:[a;+\infty)\to\mathbb{R}$ is differentiable in the usual sense, then $f'^+(t)=f'(t)$, $t\geqslant a$, where $f'(t)$ is the usual derivative of the function $f(t)$.

Consider now, for any continuous functional $V(t,\varphi)$ and solution
$x(\tau,t,\varphi)$ of equation $(1)$ with any initial conditions $(t,\varphi)$, the continuous function $v(\tau)=V(\tau,x_\tau(t,\varphi))$, $\tau\geqslant t$. The derivative of the continuous functional $V:\mathbb{R}^+\times C([-h;0],\mathbb{R}^n)\to\mathbb{R}$
along equation $(1)$ is defined as the mapping $\dot{V}:\mathbb{R}^+\times C([-h;0],\mathbb{R}^n)\to\mathbb{R}\cup\{\pm\infty\}$ given by the following formula:
$$\dot{V}(t,\varphi)=v'^+(\tau)|_{\tau=t}=v'^+(t)=
\overline{\lim_{\delta\to\, 0^+}}\,\,\frac{1}{\delta}\left\{V(t+\delta,x_{t+\delta}(t,\varphi))-V(t,\varphi)\right\}.$$

Let us now proceed to the consideration of various types of unbounded solutions of equation $(1)$.

First, we recall (see, for example, \cite{Io1}) that a solution $x(t,t_0,\varphi_0)$ of equation $(1)$ is called bounded if for this solution
there exists a number $\beta>0$ such that for all $t\geqslant t_0$ the condition $\|x_t(t_0,\varphi_0)\|\leqslant\beta$ is satisfied or, equivalently,
for all $t\geqslant t_0-h$ the condition $\|x(t,t_0,\varphi_0)\|\leqslant\beta$ is satisfied.

{\bf Definition 1}. A solution $x(t,t_0,\varphi_0)$ of equation $(1)$ will be called unbounded if this solution is not bounded, i.e., if for any number $\beta>0$ there exists a number $t_1\geqslant t_0$ such that the condition $\|x_{t_1}(t_0,\varphi_0)\|>\beta$ is satisfied or, equivalently, for all $t_1-h\leqslant t\leqslant t_1$ the condition $\|x(t,t_0,\varphi_0)\|>\beta$ is satisfied.

In what follows, for every element $x=(x_1,\dots,x_n)\in\mathbb{R}^n$, $n\geqslant 1$, and any
integer $1\leqslant k\leqslant n$
we shall use the notation $y=(x_1,\dots,x_k)\in\mathbb{R}^k$. Moreover, in what follows, for every $\varphi=(\varphi_1,\dots,\varphi_n)\in C([-h;0],\mathbb{R}^n)$ and any integer $1\leqslant k\leqslant n$ we shall use the notation $\varphi^{(k)}=(\varphi_1,\dots,\varphi_k)\in C([-h;0],\mathbb{R}^k)$.

{\bf Definition 2}. A solution $x(t,t_0,\varphi_0)$ of equation $(1)$ will be called bounded with respect to a part of the variables
$y=(x_1,\dots,x_k)$, $1\leqslant k\leqslant n$, or, more briefly, $y$-bounded, if for this solution
there exists a number $\beta>0$ such that for all $t\geqslant t_0$ the condition $\|y_t(t_0,\varphi_0)\|\leqslant\beta$ is satisfied, where $y_t(t_0,\varphi_0)=x_t(t_0,\varphi_0)^{(k)}$, or, equivalently,
for all $t\geqslant t_0-h$ the condition $\|y(t,t_0,\varphi_0)\|\leqslant\beta$ is satisfied.

It is clear that if in the definition of a $y$-bounded solution of equation $(1)$ we set $k=n$, then we obtain the definition of a bounded solution
of this equation.

{\bf Definition 3}. A solution $x(t,t_0,\varphi_0)$ of equation $(1)$ will be called unbounded with respect to a part of the variables
$y=(x_1,\dots,x_k)$, $1\leqslant k\leqslant n$, or, more briefly, $y$-unbounded, if this solution is not $y$-bounded, i.e., if for any number $\beta>0$ there exists a number $t_1\geqslant t_0$ such that the condition $\|y_{t_1}(t_0,\varphi_0)\|>\beta$ is satisfied, where $y_{t_1}(t_0,\varphi_0)=x_{t_1}(t_0,\varphi_0)^{(k)}$, or, equivalently, for all $t_1-h\leqslant t\leqslant t_1$ the condition $\|y(t,t_0,\varphi_0)\|>\beta$ is satisfied.

It is clear that if in the definition of a $y$-unbounded solution of equation $(1)$ we set $k=n$, then we obtain the definition of an unbounded solution
of this equation.

In what follows, a continuous functional $V:\mathbb{R}^+\times C([-h;0],\mathbb{R}^n)\to\mathbb{R}^+$ will be called bounded on a subset
$L\subseteq \mathbb{R}^+\times C([-h;0],\mathbb{R}^n)$ if there exists a number $M>0$ such that for all $(t,\varphi)\in L$ the condition $V(t,\varphi)\leqslant M$ is satisfied. Moreover, in what follows we shall denote by $a(r)$ and $c(r)$, where $r\in\mathbb{R}^+$,
arbitrary functions which satisfy, respectively, the following
conditions:

1) $a(r)\geqslant 0$ is an increasing function for which $a(0)=0$;

2) $c(r)\geqslant 0$ is a continuous non-decreasing function which is positive for $r>0$.

The fact that for the above-mentioned functions $a(r)$ and $c(r)$ conditions 1) and 2) are satisfied, respectively, will be written, for brevity, as
$a(r)\in{\cal K}$ and $c(r)\in{\cal K}$.

Let us formulate and prove, in terms of Lyapunov-Krasovsky functionals, a sufficient
condition for the existence of $y$-unbounded solutions of equation $(1)$.

{\bf Theorem 1}. Suppose that for equation $(1)$ and a fixed integer $1\leqslant k\leqslant n$ there exist an unbounded open subset
$G\subset C([-h;0],\mathbb{R}^n)$, a continuous functional $V:\mathbb{R}^+\times C([-h;0],\mathbb{R}^n)\to\mathbb{R}$ (Lyapunov-Krasovsky functional) with continuous derivative $\dot{V}:\mathbb{R}^+\times C([-h;0],\mathbb{R}^n)\to\mathbb{R}$ along equation $(1)$, a function $c(r)\in{\cal K}$ and a number $R_0\in\mathbb{R}^+$ for which the following conditions are satisfied:

1) $V(t,\varphi)>R_0$ for all $(t,\varphi)\in\mathbb{R}^+\times G$;

2) on every subset
$L^k_n(\beta)=\mathbb{R}^+\times \left(G\cap\left(B^k(\beta)\times C([-h;0],\mathbb{R}^{n-k})\right)\right)$ of the set $\mathbb{R}^+\times G$,
where $\beta>0$ and $B^k(\beta)=\{\varphi\in C([-h;0],\mathbb{R}^k)\,|\,\|\varphi\|\leqslant\beta\}$, the functional $V(t,\varphi)$, if $L^k_n(\beta)\not=\varnothing$, is bounded;

3) $\dot{V}(t,\varphi)\geqslant c(V(t,\varphi))$ for all $(t,\varphi)\in\mathbb{R}^+\times G$;

4) $V(t,\varphi)=R_0$ for all $(t,\varphi)\in\mathbb{R}^+\times \partial G$, where $\partial G$ is the boundary of the open
set $G$.

Then for any point $(t_0,\varphi_0)\in\mathbb{R}^+\times G$ the solution $x(t,t_0,\varphi_0)$ of equation $(1)$ is $y$-unbounded.

{\bf Proof}. Let us first show that for any solution $x(t,t_0,\varphi_0)$ of equation $(1)$, where $(t_0,\varphi_0)\in\mathbb{R}^+\times G$,
the trajectory $x_t(t_0,\varphi_0)$ corresponding to it in the space $C([-h;0],\mathbb{R}^n)$ is contained for all $t\geqslant t_0$ in the set $G$. From condition 1) of the theorem we obtain that for the point $(t_0,\varphi_0)=x_{t_0}(t_0,\varphi_0)\in\mathbb{R}^+\times G$ the inequality $V(t_0,\varphi_0)>R_0$ holds. Suppose now, for contradiction, that the trajectory $x_t(t_0,\varphi_0)$, $t\geqslant t_0$, leaves the set $G$. It is clear that the trajectory
$x_t(t_0,\varphi_0)$, $t\geqslant t_0$, can leave the set $G$ only by passing through some
boundary point of this set. Let $\varphi_1=x_{t_1}(t_0,\varphi_0)\in\partial G$ be the boundary point, where $t_1>t_0$, through which the trajectory $x_t(t_0,\varphi_0)$ passes, leaving the set $G$ for the first time. By condition 4) of the theorem we have $V(t_1,\varphi_1)=R_0$,
since $\varphi_1\in\partial G$. Since $x_t(t_0,\varphi_0)\in G$ for $t_0\leqslant t<t_1$, from condition 3)
of the theorem it follows that the inequality $V(t_0,\varphi_0)\leqslant V(t,x_t(t_0,\varphi_0))$ holds, where $t_0\leqslant t< t_1$. Passing in this inequality to the limit as $t\to t_1$, where $t_0\leqslant t< t_1$,
we obtain $R_0<V(t_0,\varphi_0)\leqslant V(t_1,\varphi_1)=R_0$, which is impossible. Thus,
the assumption that the trajectory $x_t(t_0,\varphi_0)$ leaves the set $G$ led to a contradiction and, consequently,
$x_t(t_0,\varphi_0)\in G$ for all $t\geqslant t_0$. Suppose now, for contradiction, that under the conditions of the theorem
the solution $x(t,t_0,\varphi_0)$ of equation $(1)$, where $(t_0,\varphi_0)\in\mathbb{R}^+\times G$, is $y$-bounded, i.e., for
this solution there exists a number $\beta>0$ such that $\|y_t(t_0,\varphi_0)\|\leqslant\beta$, i.e.,
$x_t(t_0,\varphi_0)\in B^k(\beta)\times C([-h;0],\mathbb{R}^{n-k})$, for all $t\geqslant t_0$. Since $x_t(t_0,\varphi_0)\in G$ for
all $t\geqslant t_0$, we have $x_t(t_0,\varphi_0)\in G\cap(B^k(\beta)\times C([-h;0],\mathbb{R}^{n-k}))$ for all $t\geqslant t_0$ and,
consequently, from condition 2) of the theorem we have the inequality $V(t,x_t(t_0,\varphi_0))\leqslant M$
for some number $M>0$ and any $t\geqslant t_0$. On the other hand, since by the conditions of the theorem the derivative $\dot{V}(t,\varphi)$ along system $(1)$ is a continuous functional, the function $\dot{V}(t,x_t(t_0,\varphi_0))$ of the variable $t$ is continuous and, consequently, using condition 3) of the theorem, we obtain for all $t\geqslant t_0$ the following inequality:
$$V(t,x_t(t_0,\varphi_0))=V(t_0,\varphi_0)+$$
$$+\int_{t_0}^t\dot{V}(\tau,x_\tau(t_0,\varphi_0))d\tau\geqslant V(t_0,\varphi_0)+c(V(t_0,\varphi_0))(t-t_0).$$
From this inequality and the inequality $M\geqslant V(t,x_t(t_0,\varphi_0))$
we obtain the inequality $M\geqslant V(t_0,\varphi_0)+c(V(t_0,\varphi_0))(t-t_0)$, which is valid for any
$t\geqslant t_0$, which is impossible for all $t\geqslant t_0$, since $c(V(t_0,\varphi_0))>0$. Thus, the contradiction assumption made above is false and, consequently, any solution $x(t,t_0,\varphi_0)$ of equation $(1)$, where $(t_0,\varphi_0)\in\mathbb{R}^+\times G$, is
$y$-unbounded.~~~$\blacksquare$

If in the conditions of Theorem 1 we set $k=n$, then this theorem becomes the following sufficient condition for the existence of
unbounded solutions of equation $(1)$ in terms of Lyapunov-Krasovsky functionals.

{\bf Corollary 1}. Suppose that for equation $(1)$ all the conditions of Theorem 1 are satisfied, in which the integer $k$ is replaced by the integer $n$, and condition
2) is replaced by the following condition:

2$'$) on every subset
$L^n(\beta)=\mathbb{R}^+\times \left(G\cap B^n(\beta)\right)\subset\mathbb{R}^+\times G$,
where $\beta>0$ and $B^n(\beta)=\{\varphi\in C([-h;0],\mathbb{R}^n)\,|\,\|\varphi\|\leqslant\beta\}$, the functional $V(t,\varphi)$, if
$L^n(\beta)\not=\varnothing$, is bounded.

Then for any point $(t_0,\varphi_0)\in\mathbb{R}^+\times G$ the solution $x(t,t_0,\varphi_0)$ of equation $(1)$ is
unbounded.~~~$\blacksquare$

Let us formulate and prove, in terms of Lyapunov-Krasovsky functionals,
one more sufficient condition for the existence of $y$-unbounded
solutions of equation $(1)$.

{\bf Theorem 2}. Suppose that for equation $(1)$ and a fixed integer $1\leqslant k\leqslant n$ there exist an unbounded open subset
$G\subset C([-h;0],\mathbb{R}^n)$, a continuous functional $V:\mathbb{R}^+\times C([-h;0],\mathbb{R}^n)\to\mathbb{R}$ (Lyapunov-Krasovsky functional) with continuous derivative $\dot{V}:\mathbb{R}^+\times C([-h;0],\mathbb{R}^n)\to\mathbb{R}$ along equation $(1)$, functions $a(r)\in{\cal K}$, $c(r)\in{\cal K}$ and a number $R_0\in\mathbb{R}^+$ for which the following conditions are satisfied:

1) $R_0<V(t,\varphi)\leqslant a(\|\varphi^{(k)}\|)$ for all $(t,\varphi)\in\mathbb{R}^+\times G$;

2) $\dot{V}(t,\varphi)\geqslant c(\|\varphi^{(k)}\|)$ for all $(t,\varphi)\in\mathbb{R}^+\times G$;

3) $V(t,\varphi)=R_0$ for all $(t,\varphi)\in\mathbb{R}^+\times \partial G$, where $\partial G$ is the boundary of the open
set $G$.

Then for any point $(t_0,\varphi_0)\in\mathbb{R}^+\times G$ the solution $x(t,t_0,\varphi_0)$ of equation $(1)$ is $y$-unbounded.

{\bf Proof}. Let us first note that from conditions 1) and 3) of the theorem being proved it follows that conditions 1) and 4) of Theorem 1 are satisfied. Let us now show
that from conditions 1) and 2) of the theorem being proved it follows that condition 3) of Theorem 1 is satisfied. Indeed, since $a(r)\in{\cal K}$, the
function $a(r)$ is increasing and, consequently, for this function there is an inverse function $a^{-1}(r)$, which is also increasing.
Computing the function $a^{-1}(r)$ from the left and right sides of the inequality $V(t,\varphi)\leqslant a(\|\varphi^{(k)}\|)$, we obtain the inequality $\|\varphi^{(k)}\|\geqslant a^{-1}(V(t,\varphi))$. If now we compute the function $c(r)$ from the left and right sides of this inequality, then
we obtain on the set $\mathbb{R}^+\times G$ the inequality $c(\|\varphi^{(k)}\|)\geqslant (c\circ a^{-1})(V(t,\varphi))$. From
this, using condition 2) of the theorem being proved, we obtain the required inequality $\dot{V}(t,\varphi)\geqslant\tilde{c}(V(t,\varphi))$,
where $\tilde{c}(r)=(c\circ a^{-1})(r)\in{\cal K}$. Thus, from the conditions of the theorem being proved it follows that conditions 1), 3)
and 4) of Theorem 1 are satisfied. Further, repeating verbatim the arguments carried out in the proof of Theorem 1 using conditions 1), 3) and 4), we obtain that for any point $(t_0,\varphi_0)\in\mathbb{R}^+\times G$ the solution $x(t,t_0,\varphi_0)$ of equation $(1)$ satisfies the condition
$x(t,t_0,\varphi_0)\in G$ for all $t\geqslant t_0$.
Let us now show that for any point $(t_0,\varphi_0)\in\mathbb{R}^+\times G$ the solution $x(t,t_0,\varphi_0)$ of equation $(1)$ is $y$-unbounded. Suppose, for contradiction, that under the conditions of the theorem
the solution $x(t,t_0,\varphi_0)$ of equation $(1)$, where $(t_0,\varphi_0)\in\mathbb{R}^+\times G$, is $y$-bounded, i.e., for
this solution there exists a number $\beta>0$ such that for all $t\geqslant t_0$ the condition $\|y_t(t_0,\varphi_0)\|\leqslant\beta$ is satisfied.
Using condition 1) of the theorem being proved, condition 3) of Theorem 1 and the condition $a(r)\in{\cal K}$, analogously to what was done in the proof of Theorem 1, we obtain for
the above-mentioned function $\tilde{c}(r)$ and all $t\geqslant t_0$ the inequality
$a(\beta)\geqslant V(t_0,\varphi_0)+\tilde{c}(V(t_0,\varphi_0))(t-t_0)$, which is impossible for all $t\geqslant t_0$, since $\tilde{c}(V(t_0,\varphi_0))>0$. Thus, the contradiction assumption made above is false and, consequently, any solution $x(t,t_0,\varphi_0)$ of equation $(1)$, where $(t_0,\varphi_0)\in\mathbb{R}^+\times G$, is
$y$-unbounded.~~~$\blacksquare$

If in the conditions of Theorem 2 we set $k=n$, then this theorem becomes the following sufficient condition for the existence of
unbounded solutions $x(t,t_0,\varphi_0)$ of equation $(1)$
in terms of Lyapunov-Krasovsky functionals.

{\bf Corollary 2}. Suppose that for equation $(1)$ all the conditions
of Theorem 2 are satisfied, in which the integer $k$ is replaced by the integer $n$, and conditions 1) and
2) are replaced by the following conditions:

1$'$) $R_0<V(t,\varphi)\leqslant a(\|\varphi\|)$ for all $(t,\varphi)\in\mathbb{R}^+\times G$;

2$'$) $\dot{V}(t,\varphi)\geqslant c(\|\varphi\|)$ for all $(t,\varphi)\in\mathbb{R}^+\times G$.

Then for any point $(t_0,\varphi_0)\in\mathbb{R}^+\times G$ the solution $x(t,t_0,\varphi_0)$ of equation $(1)$ is unbounded.~~~$\blacksquare$

Let us now obtain a sufficient condition for the existence of $y$-unbounded solutions of equation $(1)$ using two
Lyapunov-Krasovsky functionals.

{\bf Theorem 3}. Suppose that for equation $(1)$ and a fixed integer $1\leqslant k\leqslant n$ there exist an unbounded open subset
$G\subset C([-h;0],\mathbb{R}^n)$, a continuous functional $W:\mathbb{R}^+\times C([-h;0],\mathbb{R}^n)\to\mathbb{R}$ (Lyapunov-Krasovsky functional), a continuous functional $V:\mathbb{R}^+\times C([-h;0],\mathbb{R}^k)\to\mathbb{R}$ (Lyapunov-Krasovsky functional) with continuous derivative $\dot{V}:\mathbb{R}^+\times C([-h;0],\mathbb{R}^n)\to\mathbb{R}$ along equation $(1)$, a function $c(r)\in{\cal K}$ and a number $R_0\in\mathbb{R}^+$ for which the following conditions are satisfied:

1) $W(t,\varphi)>R_0$ for all $(t,\varphi)\in\mathbb{R}^+\times G$;

2) $\dot{W}(t,\varphi)\geqslant 0$ for all $(t,\varphi)\in\mathbb{R}^+\times G$;

3) $W(t,\varphi)=R_0$ for all $(t,\varphi)\in\mathbb{R}^+\times\partial G$;

4) $V(t,\varphi^{(k)})>0$ for all $(t,\varphi)\in\mathbb{R}^+\times G$;

5) On every subset $N^k(\beta)=\mathbb{R}^+\times (G_k\cap B^k(\beta))\subset\mathbb{R}^+\times G_k$, where $\beta>0$ and $G_k=\{\psi\in C([-h;0],\mathbb{R}^k)\,|\,\psi=\varphi^{(k)}, \varphi\in G\}$, the functional $V(t,\psi)$, if $N^k(\beta)\not=\varnothing$, is bounded;

6) $\dot{V}(t,\varphi)\geqslant c(V(t,\varphi^{(k)}))$ for all $(t,\varphi)\in\mathbb{R}^+\times G$.

Then for any point $(t_0,\varphi_0)\in\mathbb{R}^+\times G$ the solution $x(t,t_0,\varphi_0)$ of equation $(1)$ is $y$-unbounded.

{\bf Proof}. Analogously to what was done in the proof of Theorem 1, using the functional $W(t,\varphi)$ and conditions 1)\,--\,3)
of the theorem being proved, it is shown that for any solution $x(t,t_0,\varphi_0)$ of equation $(1)$, where $(t_0,\varphi_0)\in\mathbb{R}^+\times G$,
the trajectory $x_t(t_0,\varphi_0)$ corresponding to it in the space $C([-h;0],\mathbb{R}^n)$ is contained for all $t\geqslant t_0$ in the set $G$.
Now, analogously to what was done in the proof of Theorem 1, using the functional $V(t,\psi)$ and conditions 4)\,--\,6) of the theorem being proved, it is shown that the solution $x(t,t_0,\varphi_0)$ of equation $(1)$, where $(t_0,\varphi_0)\in\mathbb{R}^+\times G$, cannot be $y$-bounded. Thus, it is shown that for any point $(t_0,\varphi_0)\in\mathbb{R}^+\times G$ the solution $x(t,t_0,\varphi_0)$ of equation $(1)$ is $y$-unbounded.~~~$\blacksquare$

If in the conditions of Theorem 3 we set $k=n$, then this theorem becomes the following sufficient condition for the existence of
unbounded solutions of equation $(1)$ in terms of two Lyapunov-Krasovsky functionals.

{\bf Corollary 3}. Suppose that for equation $(1)$ all the conditions
of Theorem 3 are satisfied, in which the integer $k$ is replaced by the integer $n$, and conditions 4)\,--\,6) are replaced, respectively, by the following
conditions:

4$'$) $V(t,\varphi)>0$ for all $(t,\varphi)\in\mathbb{R}^+\times G$;

5$'$) On every subset $N^n(\beta)=\mathbb{R}^+\times (G\cap B^n(\beta))\subset\mathbb{R}^+\times G$, where $\beta>0$, the functional
$V(t,\varphi)$, if $N^n(\beta)\not=\varnothing$, is bounded;

6$'$) $\dot{V}(t,\varphi)\geqslant c(V(t,\varphi))$ for all $(t,\varphi)\in\mathbb{R}^+\times G$.

Then for any point $(t_0,\varphi_0)\in\mathbb{R}^+\times G$ the solution $x(t,t_0,\varphi_0)$
of equation $(1)$ is unbounded.~~~$\blacksquare$

Let us formulate and prove, in terms of two Lyapunov-Krasovsky functionals,
one more sufficient condition for the existence of $y$-unbounded
solutions $x(t,t_0,\varphi_0)$ of equation $(1)$.

{\bf Theorem 4}. Suppose that for equation $(1)$ and a fixed integer $1\leqslant k\leqslant n$ there exist an unbounded open subset
$G\subset C([-h;0],\mathbb{R}^n)$, a continuous functional $W:\mathbb{R}^+\times C([-h;0],\mathbb{R}^n)\to\mathbb{R}$ (Lyapunov-Krasovsky functional) and a number $R_0\in\mathbb{R}^+$ satisfying conditions 1)\,--\,3) of Theorem 3. Moreover, suppose that there exist a continuous functional $V:\mathbb{R}^+\times C([-h;0],\mathbb{R}^k)\to\mathbb{R}$ (Lyapunov-Krasovsky functional) with continuous derivative $\dot{V}:\mathbb{R}^+\times C([-h;0],\mathbb{R}^n)\to\mathbb{R}$ along equation $(1)$ and functions $a(r),c(r)\in{\cal K}$ for which the following conditions are satisfied:

1) $0<V(t,\varphi^{(k)})\leqslant a(\|\varphi^{(k)}\|)$ for all $(t,\varphi)\in\mathbb{R}^+\times G$;

2) $\dot{V}(t,\varphi)\geqslant c(\|\varphi^{(k)}\|)$ for all $(t,\varphi)\in\mathbb{R}^+\times G$.

Then for any point $(t_0,\varphi_0)\in\mathbb{R}^+\times G$ the solution $x(t,t_0,\varphi_0)$ of equation $(1)$ is $y$-unbounded.

{\bf Proof}. Since the continuous functional $W(t,\varphi)$ satisfies conditions 1)\,--\,3)
of Theorem 3, analogously to what was done in the proof of
Theorem 3, we obtain that for any solution $x(t,t_0,\varphi_0)$ of equation $(1)$, where $(t_0,\varphi_0)\in\mathbb{R}^+\times G$,
the trajectory $x_t(t_0,\varphi_0)$ corresponding to it in the space $C([-h;0],\mathbb{R}^n)$ is contained for all $t\geqslant t_0$ in the set $G$.
Now, analogously to what was done in the proof of Theorem 2, using conditions 1) and 2) of the theorem being proved, it is shown that
for any point $(t_0,\varphi_0)\in\mathbb{R}^+\times G$ the solution $x(t,t_0,\varphi_0)$ of equation $(1)$ is $y$-unbounded.~~~$\blacksquare$

If in the conditions of Theorem 4 we set $k=n$, then this theorem becomes the following sufficient condition for the existence of
unbounded solutions $x(t,t_0,\varphi_0)$ of equation $(1)$
in terms of two Lyapunov-Krasovsky functionals.

{\bf Corollary 4}. Suppose that for equation $(1)$ all the conditions
of Theorem 4 are satisfied, in which the integer $k$ is replaced by the integer $n$, and conditions 1) and 2) are replaced, respectively, by the following
conditions:

1$'$) $0<V(t,\varphi)\leqslant a(\|\varphi\|)$ for all $(t,\varphi)\in\mathbb{R}^+\times G$;

2$'$) $\dot{V}(t,\varphi)\geqslant c(\|\varphi\|)$ for all $(t,\varphi)\in\mathbb{R}^+\times G$.

Then for any point $(t_0,\varphi_0)\in\mathbb{R}^+\times G$ the solution $x(t,t_0,\varphi_0)$ of equation $(1)$ is unbounded.~~~$\blacksquare$

The following statement is a refinement of Theorem 2 for the case when some additional requirements connected with the continuability of the $y$-components of these solutions are imposed on the solutions of equation $(1)$.

{\bf Proposition 1}. Suppose that for equation $(1)$ all the conditions of Theorem 2 are satisfied and, moreover, suppose that for every solution $x(t,t_0,\varphi_0)$, where $(t_0,\varphi_0)\in\mathbb{R}^+\times G$, of equation $(1)$ the mapping $y_t(t_0,\varphi_0)$ corresponding to it is continuable to the whole half-line $[t_0;+\infty)$. Then for any point $(t_0,\varphi_0)\in\mathbb{R}^+\times G$ the solution $x(t,t_0,\varphi_0)$ of equation $(1)$ satisfies the condition
$\lim\limits_{t\to+\infty}\|y_t(t_0,\varphi_0)\|=+\infty$.

{\bf Proof}. Analogously to what was done in the proofs of Theorems 1 and 2, for any solution
$x(t,t_0,\varphi_0)$, where $(t_0,\varphi_0)\in\mathbb{R}^+\times G$, of equation $(1)$ and the function $\tilde{c}(r)$ defined in the proof of Theorem 2, we obtain for all $t\in[t_0;+\infty)$ the inequality
$a(\|y_t(t_0,\varphi_0)\|)\geqslant  V(t_0,\varphi_0)+\tilde{c}(V(t_0,\varphi_0))(t-t_0)$.
Since $\tilde{c}(V(t_0,\varphi_0))>0$, we have $\lim\limits_{t\to+\infty}(V(t_0,\varphi_0)+\tilde{c}(V(t_0,\varphi_0))(t-t_0))=+\infty$
and, consequently, $\lim\limits_{t\to+\infty}a(\|y_t(t_0,\varphi_0)\|)=+\infty$.
From this we obtain $\lim\limits_{t\to+\infty}\|y_t(t_0,\varphi_0)\|=+\infty$, since $a(r)\in{\cal K}$. Thus, it is shown that
under the requirements of the statement being proved, for any point $(t_0,\varphi_0)\in\mathbb{R}^+\times G$ the solution $x(t,t_0,\varphi_0)$ of equation $(1)$ satisfies the condition
$\lim\limits_{t\to+\infty}\|y_t(t_0,\varphi_0)\|=+\infty$.~~~$\blacksquare$

If in the conditions of Proposition 1 we set $k=n$, then this statement becomes the following refinement of Corollary 2 for the case when some additional requirements connected with the continuability of these solutions are imposed on the solutions of equation $(1)$.

{\bf Corollary 5}. Suppose that for equation $(1)$ all the conditions of Corollary 2 are satisfied and, moreover, suppose that for every solution $x(t,t_0,\varphi_0)$, where $(t_0,\varphi_0)\in\mathbb{R}^+\times G$, of equation $(1)$ the mapping $x_t(t_0,\varphi_0)$ corresponding to it is continuable to the whole half-line $[t_0;+\infty)$. Then for any point $(t_0,\varphi_0)\in\mathbb{R}^+\times G$ the solution $x(t,t_0,\varphi_0)$ of equation $(1)$ satisfies the condition
$\lim\limits_{t\to+\infty}\|x_t(t_0,\varphi_0)\|=+\infty$.~~~$\blacksquare$

The following statement, whose proof is analogous to the proof of Proposition 1, is a refinement of Theorem 4 for the case when some additional requirements connected with the continuability of the $y$-components of these solutions are imposed on the solutions of equation $(1)$.

{\bf Proposition 2}. Suppose that for equation $(1)$ all the requirements of Theorem 4 are satisfied and, moreover, suppose that for every solution $x(t,t_0,\varphi_0)$, where $(t_0,\varphi_0)\in\mathbb{R}^+\times G$, of equation $(1)$ the mapping $y_t(t_0,\varphi_0)$ corresponding to it is continuable to the whole half-line $[t_0;+\infty)$. Then for any point $(t_0,\varphi_0)\in\mathbb{R}^+\times G$ the solution $x(t,t_0,\varphi_0)$ of equation $(1)$ satisfies the condition
$\lim\limits_{t\to+\infty}\|y_t(t_0,\varphi_0)\|=+\infty$.~~~$\blacksquare$

If in the conditions of Proposition 2 we set $k=n$, then this statement becomes the following refinement of Corollary 4 for the case when some additional requirements connected with the continuability of these solutions are imposed on the solutions of equation $(1)$.

{\bf Corollary 6}. Suppose that for equation $(1)$ all the conditions of Corollary 4 are satisfied and, moreover, suppose that for every solution $x(t,t_0,\varphi_0)$, where $(t_0,\varphi_0)\in\mathbb{R}^+\times G$, of equation $(1)$ the mapping $x_t(t_0,\varphi_0)$ corresponding to it is continuable to the whole half-line $[t_0;+\infty)$. Then for any point $(t_0,\varphi_0)\in\mathbb{R}^+\times G$ the solution $x(t,t_0,\varphi_0)$ of equation $(1)$ satisfies the condition
$\lim\limits_{t\to+\infty}\|x_t(t_0,\varphi_0)\|=+\infty$.~~~$\blacksquare$

In conclusion, let us give an example of the application of Theorem 1. Let an integro-differential system be given
$$\left\{\begin{array}{l}
\dot{x}_1(t)=e^hx(t-h)+\sin^2(x_2(t-h))\left(x_1(t)+\int_{t-h}^{t}e^{t-s}x(s)\,ds\right)^2\\
\dot{x}_2(t)=x_1(t)+x_1^2(t-h)+x_2(t-h)\\
\end{array}\right.\!\!\!\!,~t\in\mathbb{R}^+,$$
with constant delay $h>0$. Let us define for this system the Lyapunov-Krasovsky functional
$V:\mathbb{R}^+\times C([-h;0],\mathbb{R}^2)\to\mathbb{R}$, by setting
$$V(t,\varphi)=V_1(\varphi_1)=\left(\varphi_1(0)+\int_{-h}^{0}e^{-\vartheta}\varphi_1(\vartheta)\,d\vartheta\right)^2,
~\mbox{where}$$
$$\varphi=(\varphi_1,\varphi_2)\in C([-h,0],\mathbb{R})\times C([-h,0],\mathbb{R})=C([-h,0],\mathbb{R}^2),$$
and the infinite open subset $G\subset C([-h,0],\mathbb{R}^2)$, by setting
$$G=\{\varphi\in C([-h,0],\mathbb{R}^2)\mid V(t,\varphi)>0\}=G_1\times C([-h;0],\mathbb{R}),$$
where $G_1=\{\varphi_1\in C([-h;0],\mathbb{R})\mid V_1(\varphi_1)>0\}$.
Let us show, using Theorem 1, that for any point $(t_0,\varphi_0)\in\mathbb{R}^+\times G$ the solution $x(t,t_0,\varphi_0)$ of the given system is $y$-unbounded, where $y=(x_1)$. Let us set $R_0=0$. Then from the definition of the set $G$ it follows that $V(t,\varphi)>R_0$ for all $(t,\varphi)\in G$.
Thus, condition 1) of Theorem 1 is satisfied. It is easy to see that for all $(t,\varphi)\in\partial G$ the equality $V(t,\varphi)=R_0$ is valid, since the boundary $\partial G$ of the set $G$ is the set $\partial G=\{\varphi\in C([-h,0],\mathbb{R}^2)\mid V(t,\varphi)=0\}$. Thus, condition 4) of Theorem 1 is satisfied. Let us verify the fulfillment of condition 2) of Theorem 1. Consider, for any $\beta>0$, the subset
$$L^1_2(\beta)=\mathbb{R}^+\times \left(G\cap\left(B^1(\beta)\times C([-h;0],\mathbb{R}^{2-1})\right)\right)=$$
$$=\mathbb{R}^+\times \left(G_1\cap B^1(\beta)\right)\times C([-h;0],\mathbb{R})\subset\mathbb{R}^+\times G,$$
where $B^1(\beta)=\{\varphi_1\in C([-h;0],\mathbb{R})\mid \|\varphi_1\|\leqslant\beta\}$.
Since the boundary $\partial G_1$ of the set $G_1$ is the set
$\partial G_1=\{\varphi_1\in C([-h;0],\mathbb{R})\mid V_1(\varphi_1)=0\}$,
and since the function $\varphi_1=0$ belongs to $\partial G_1\cap B^1(\beta)$, we have $G_1\cap B^1(\beta)\not=\varnothing$ and, consequently, $L^1_2(\beta)\not=\varnothing$.
Let us show that the functional $V(t,\varphi)$ is bounded on the set $L^1_2(\beta)$. Indeed, for any element $(t,\varphi)\in L^1_2(\beta)$ we have
$$V(t,\varphi)=V_1(\varphi_1)=\left(\varphi_1(0)+\int_{-h}^{0}e^{-\vartheta}\varphi_1(\vartheta)\,d\vartheta\right)^2=
\left|\varphi_1(0)+\int_{-h}^{0}e^{-\vartheta}\varphi_1(\vartheta)\,d\vartheta\right|^2\leqslant$$
$$\leqslant \left(|\varphi_1(0)|+\int_{-h}^{0}e^{-\vartheta}|\varphi_1(\vartheta)|\,d\vartheta\right)^2
\leqslant \left(\|\varphi_1\|+\int_{-h}^{0}e^{-\vartheta}\|\varphi_1\|\,d\vartheta\right)^2\leqslant$$
$$\leqslant \left(\beta+\beta\int_{-h}^{0}e^{-\vartheta}\,d\vartheta\right)^2=(\beta(1+(-1+e^h))^2=\beta^2e^{2h}=M>0.$$
Thus, the functional $V(t,\varphi)$ is bounded on the set $L^1_2(\beta)$ and, consequently, condition 2 of Theorem 1 is satisfied.
Let us verify the fulfillment of condition 3) of Theorem 1. Consider, for an arbitrary point $(t,\varphi)\in\mathbb{R}^+\times C([-h;0],\mathbb{R}^2)$,
the derivative $\dot{V}(t,\varphi)=v'(\tau)|_{\tau=t}=v'(t)$ along the given system, where
$$v(\tau)=V(\tau,x_\tau)=\left(x_1(\tau)+\int_{-h}^{0}e^{-\vartheta}x_1(\tau+\vartheta)\,d\vartheta\right)^2=$$
$$=\left(x_1(\tau)+\int_{-h}^{0}e^{-\vartheta}x_1(\tau+\vartheta)\,d\vartheta\right)^2=
\left(x_1(\tau)+\int_{\tau-h}^{\tau}e^{-(s-\tau)}x_1(s)\,ds\right)^2$$
and $x(\tau,t,\varphi)=(x_1(\tau),x_2(\tau))$ is the solution of the given system with initial conditions $(t,\varphi)$. Since
$v(t)=\xi^2(t)$, where
$\xi(t)=x_1(t)+e^t\!\int_{t-h}^{t}e^{-s}x_1(s)\,ds$, then
for the derivative $v'(t)$ we have
$$v'(t)=2\xi(t)\xi'(t)=2\xi(t)(e^hx_1(t-h)+\xi(t)\sin^2(x_2(t-h))+e^t\!\int_{t-h}^{t}e^{-s}x_1(s)\,ds\,+$$
$$+\,x_1(t)-e^hx_1(t-h))=2v(t)(\sin^2(x_2(t-h))+1)\geqslant 2v(t).$$
Since $v'(t)=\dot{V}(t,\varphi)$ and $v(t)=v(\tau)|_{\tau=t}=V(\tau,x_\tau)|_{\tau=t}=V(t,\varphi)$, for any point $(t,\varphi)\in\mathbb{R}^+\times C([-h;0],\mathbb{R}^2)$ and, in particular, for all $(t,\varphi)\in\mathbb{R}^+\times G$
we have the inequality $\dot{V}(t,\varphi)\geqslant 2V(t,\varphi)$. Therefore, if as the function
$c(r)\in{\cal K}$ we take the function $c(r)=2r$, then we obtain $\dot{V}(t,\varphi)\geqslant c(V(t,\varphi))$ and, consequently, condition 3) of Theorem 1 is satisfied. Thus, for the given system all the conditions of Theorem 1 are satisfied. It follows that for any point $(t_0,\varphi_0)\in\mathbb{R}^+\times G$ the solution $x(t,t_0,\varphi_0)$ of the given system is $y$-unbounded, where $y=(x_1)$.~~~$\blacktriangledown$

\vspace{2cm}

\noindent Lapin K.S. - Mordovian State Pedagogical University named after M.E. Evseviev, Saransk, Russia,
e-mail: klapin@mail.ru
\vspace{0.5cm}

\noindent Lapin S.V. - Serov Str., Saransk, Russia,
e-mail: slapin@mail.ru

\end{document}